\documentclass[a4paper,12pt]{amsart}
\usepackage{amssymb}
\usepackage{ifthen}
\nonstopmode \numberwithin{equation}{section}
\newtheorem{thm}{Theorem}
\newtheorem{cor}{Corollary}
\newtheorem{lem}{Lemma}

\newtheorem{conj}{Conjecture}

\theoremstyle{definition}
\newtheorem{defn}{Definition}[section]
\newtheorem{example}{Example}[section]
\newtheorem{prob}[equation]{Problem}

\newenvironment{rem}{%
\bigskip
\noindent \textsl{{\sl Remark. }}}{\bigskip}
\newenvironment{rems}{%
\bigskip
\noindent \textsl{{\sl Remarks. }}}{\bigskip}

\newcounter {own}
\def\theown {\thesection       .\arabic{own}}

\newenvironment{pf}[1][]{%
 \vskip 3mm
 \noindent
 \ifthenelse{\equal{#1}{}}%
  {{\slshape Proof. }}%
  {{\slshape #1.} }%
 }%
{\qed\bigskip}

\newcounter{alphabet}
\newcounter{tmp}
\newenvironment{Thm}[1][]{\refstepcounter{alphabet}%
\bigskip%
\noindent%
{\bf Theorem \Alph{alphabet}}%
\ifthenelse{\equal{#1}{}}{}{ (#1)}%
{\bf .} \itshape}{\vskip 8pt}

\newcommand{\ID}{{\mathbb D}}
\newcommand{\IDs}{{\mathbb D^*}}
\newcommand{\IDb}{{\overline{\mathbb D}}}
\newcommand{\IDsb}{{\overline{\mathbb D^*}}}

\newcommand{\IC}{{\mathbb C}}

\newcommand{\sphere}{{\widehat{\mathbb C}}}

\newcommand{\inv}{^{-1}}

\newcommand{\arctanh}{{\operatorname{arctanh}}}

\def\be{\begin{equation}}
\def\ee{\end{equation}}

\newcommand{\bee}{\begin{enumerate}}
\newcommand{\eee}{\end{enumerate}}

\newcommand{\blem}{\begin{lem}}
\newcommand{\elem}{\end{lem}}
\newcommand{\bthm}{\begin{thm}}
\newcommand{\ethm}{\end{thm}}
\newcommand{\bcor}{\begin{cor}}
\newcommand{\ecor}{\end{cor}}
\newcommand{\beg}{\begin{examp}}
\newcommand{\eeg}{\end{examp}}
\newcommand{\begs}{\begin{examples}}
\newcommand{\eegs}{\end{examples}}
\newcommand{\bdefe}{\begin{defn}}
\newcommand{\edefe}{\end{defn}}
\newcommand{\bprob}{\begin{prob}}
\newcommand{\eprob}{\end{prob}}
\newcommand{\bei}{\begin{itemize}}
\newcommand{\eei}{\end{itemize}}

\newcommand{\bcon}{\begin{conj}}
\newcommand{\econ}{\end{conj}}
\newcommand{\bcons}{\begin{conjs}}
\newcommand{\econs}{\end{conjs}}
\newcommand{\bprop}{\begin{propo}}
\newcommand{\eprop}{\end{propo}}
\newcommand{\br}{\begin{rem}}
\newcommand{\er}{\end{rem}}
\newcommand{\brs}{\begin{rems}}
\newcommand{\ers}{\end{rems}}
\newcommand{\bo}{\begin{obser}}
\newcommand{\eo}{\end{obser}}
\newcommand{\bos}{\begin{obsers}}
\newcommand{\eos}{\end{obsers}}
\newcommand{\bpf}{\begin{pf}}
\newcommand{\epf}{\end{pf}}
\newcommand{\ba}{\begin{array}}
\newcommand{\ea}{\end{array}}
\newcommand{\beq}{\begin{eqnarray}}
\newcommand{\beqq}{\begin{eqnarray*}}
\newcommand{\eeq}{\end{eqnarray}}
\newcommand{\eeqq}{\end{eqnarray*}}

\newcounter{minutes}
\divide\time by 60
\newcounter{hours}
\multiply\time by 60 \addtocounter{minutes}{-\time}

\begin{document}
\bibliographystyle{amsplain}
\title{Loewner chain and quasiconformal extension of some classes of univalent functions}
\begin{center}
{\tiny \texttt{FILE:~\jobname .tex,
        printed: \number\year-\number\month-\number\day,
        \thehours.\ifnum\theminutes<10{0}\fi\theminutes}
}
\end{center}

\author{Bappaditya Bhowmik${}^{~\mathbf{*}}$}
\address{Bappaditya Bhowmik, Department of Mathematics,
Indian Institute of Technology Kharagpur, Kharagpur - 721302, India.}
\email{bappaditya@maths.iitkgp.ac.in}
\author{Goutam Satpati}
\address{Goutam Satpati, Department of Mathematics,
Indian Institute of Technology Kharagpur, Kharagpur - 721302, India.}
\email{gsatpati@iitkgp.ac.in}

\subjclass[2010]{30C62, 30C55}
\keywords{ Loewner Chain, Meromorphic, Quasiconformal.\\
${}^{\mathbf{*}}$ Corresponding author}
\date{ \today ; File: BS-v5-Loewner.tex}
\begin{abstract}
In this article we obtain quasiconformal extensions of some classes of conformal maps
defined either on the unit disc or on the exterior of it onto the extended complex plane.
Some of these extensions have been obtained by constructing suitable Loewner chains
and others have been obtained by applying a well-known result.
\end{abstract}

\thanks{The first author of this article would like to thank
CSIR, India (Ref.No.- 25(0281)/18/EMR-II) for its financial support.}

\maketitle
\pagestyle{myheadings}
\markboth{Loewner Chain}{Loewner Chain}
\bigskip

\section{Introduction}
In 1923 Loewner developed a theory to represent conformal maps parametrically-known as Loewner Theory. Later this
theory was developed by Kufarev and Pommerenke. It describes a family of time parameterized
conformal maps defined on the unit disc $\ID:=\{z\in \IC : |z|<1\}$
(here and hereafter $\IC$ denotes the complex plane)
whose images are simply connected domains continuously increasing with time.
Let $f_t(z)=f(z,t)= e^tz+ \sum_{n=2}^{\infty}a_n(t)z^n$ be a function defined on
$\ID \times [0,\infty)$ where $z\in \ID$ and  $t\in [0,\infty)$. Such a function
$f_t(z)$ is called a
\textit{Loewner chain} if $f_t(z)$ is analytic
and univalent in $\ID$ for each fixed $t\in [0,\infty)$ and satisfies
$f_s(\ID) \subsetneq f_t(\ID)$ for $0\leq s<t< \infty$.
In 1965, Pommerenke (see \cite[Theorem\,6.2]{PM2}) proved the following necessary
and sufficient condition for Loewner chain.

\begin{Thm}\label{thA}
Let $0<r_0\leq1$ and $\ID_{r_0}= \lbrace z: |z|< r_0\rbrace$.
The function $f(z,t)= e^tz+ \sum_{n=2}^{\infty}a_n(t)z^n$ defined on
$\ID \times [0,\infty)$ is a Loewner chain if and only if the following two
conditions are satisfied :
\begin{itemize}
\item[(i)] $f(z,t)$ is analytic in $z \in \ID_{r_0}$, for each $t \in [0,\infty)$, absolutely continuous in $t$ for each $z \in \ID_{r_0}$ and satisfies
$$
|f(z,t)| \leq K_0 e^t, \quad z\in \ID_{r_0},\,t \in [0,\infty),
$$
where $K_0$ is a positive constant.
\item[(ii)] There exists a function $p(z,t)$ analytic in $z\in \ID$, for each $t$ and
measurable in $t$ for each $z\in \ID$, satisfying
$\mathrm{Re} ~p(z,t)>0$, for all $z\in \ID$ and $t\in [0,\infty)$, such that
\begin{equation}\label{eq1}
\dot{f}(z,t)= z f'(z,t)p(z,t), \quad z\in \ID_{r_0}, ~a.e.~
t\in [0,\infty),
\end{equation}
where $\dot{f}=\partial f/\partial t$ and $f'= \partial f/ \partial z$.
\end{itemize}
\end{Thm}
The partial differential equation \eqref{eq1} is known as \textit{Loewner-Kufarev PDE} and the function $p(z,t)$ is called \textit{Herglotz function}. Here we
mention that the term $e^t=f'(0,t)$ may be replaced by any complex valued function
$a_1(t)\neq 0$ such that it is locally absolutely continuous on $[0,\infty)$
and $\lim_{t \to \infty}|a_1(t)|= \infty $ (see \cite[Theorem~A$'$]{IH2}). In this case,
the inequality in the first condition of Theorem A will be replaced by
$|f(z,t)| \leq K_0|a_1(t)|$ for all $z\in \ID_{r_0}$ and $t\in [0, \infty)$.
The Loewner chain $f(z,t)$ with $a_1(t)=e^t$ is called  \textit{standard
Loewner chain}.

The method of Loewner chain plays a central role in the de Brange's proof of famous
Bieberbach conjecture.
This method also solves many important problems in univalent function theory which was not
accessible by other existing methods.
The theory of Loewner chain has recently been  developed in
many other directions which include a new innate approach suggested by Bracci, Contreras and Díaz-
Madrigal (\cite{BCD-1},\cite{BCD-2}) and the Schramm's stochastic version of the Loewner differential
equation  (\cite{OS}).
Another important aspect of Loewner theory lies in
quasiconformal extension of conformal maps defined on $\ID$, to the Riemann Sphere
$\sphere= \IC \cup \{\infty\}$, which was found by Becker in 1972. In this article our interest lies in the later one.

A sense preserving homeomorphism $f: \sphere\to \sphere$ is called a
$k$-\textit{quasiconformal mapping} if it belongs to the Sobolev
space $W^{1,2}_{\textrm{loc}}(\sphere)$ and satisfies the inequality
$|f_{\bar z}| \leq k |f_z|$ a.e., for some $k\in [0,1)$, where $f_{\bar z}:= \partial f/\partial \overline{z}$ and $f_z:= \partial f/\partial z$. The quantity $\mu_f:=f_{\bar z}/f_z$ is called \textit{complex dilatation} of $f$. We mention here that such a function $f$ is also called $K$-quasiconformal
in the literature, where $K=(1+k)/(1-k) \geq 1$. Let $f$ be a conformal map defined
on a domain $\Omega \subset \IC$. We say that $f:\Omega \to \IC$ has a
$k$-\textit{quasiconformal extension} onto $\IC$ or onto $\sphere$, if there exists a $k$-quasiconformal mapping $F: \IC \to \IC$, or $F: \sphere \to \sphere$
respectively,
such that $F|_{\Omega}=f$.
By the removability property of quasiconformal mappings (see f.i.\,\cite[p.\,45]{Lehto0}),
it follows that a function $f:\Omega \to \IC$ is $k$-quasiconformally extendible
to $\IC$ if and only if it admits a $k$-quasiconformal extension $F$ onto $\sphere$
with $F(\infty)=\infty$. Normalized holomorphic maps defined on $\ID$ having quasiconformal extension to $\sphere$ play a very
important role in
Teichm{\" u}ller theory as they can be identified with the elements of the universal
Teichm{\" u}ller space (see e.g.\,\cite[chap.\,III]{Lehto2}). J. Becker
(compare \cite{JB1,JB2}) found a remarkable result concerning quasiconformal extension of univalent functions with the help of Loewner chain.
We first state Becker's
result in this direction (see f.i.\,\cite[Theorem 5.2]{JB2}).

\begin{Thm}\label{thB}
Suppose that $f_t(z)=f(z,t),~(z,t)\in \ID \times [0, \infty)$ be a standard Loewner chain such that the Herglotz function $p(z,t)$ in \eqref{eq1} lies within $\mathcal{D}(k)$ for all $z\in \ID$ and almost all $t\in [0,\infty)$, where $\mathcal{D}(k)$ is the closed hyperbolic disc in the right half plane centered at $1$ and radius $\arctanh \,k$, given by
$$
\mathcal{D}(k) = \left\{w\in \IC : \left|\frac{w-1}{w+1} \right| \leq k \right \}
= \left\{w\in \IC : \left| w - \frac{1+k^2}{1-k^2} \right| \leq \frac{2k}{1-k^2} \right\}.
$$
Then for each $t\geq 0$ the function $f_t$ admits a $k$-quasiconformal extension
onto $\sphere$ and such an extension for $f_0$ is given by
\begin{equation}\label{eq2}
F(z)=
\begin{cases}
f_0(z)=f(z,0), & z\in \ID \\
f(z/|z|,\log |z|) & z \in \IDsb ,
\end{cases}
\end{equation}
which fixes infinity i.e. $F(\infty)=\infty$, where
$\overline{\ID^{*}}= \{z : |z|\geq 1\}$.
\end{Thm}

If $f(z,t)$ be a Loewner chain with the complex first coefficient $a_1(t)$, then also
Theorem B is valid (see f.i. \cite[Theorem\,C$^{\prime}$]{IH2}).
Many sufficient conditions for univalence and quasiconformal extensibility have been proved using Theorem B together with the Pommerenke's criteria (Theorem A) for Loewner chain. Quasiconformal extension criteria for starlike, convex, close-to-convex, spirallike and Bazilevi{\v c} functions have been found by this technique where one can easily obtain explicit extended quasiconformal maps from the suitable Loewner chains, which is the main advantage of this method (compare \cite{JB2,IH1,IH2}).

Now we introduce some important classes of univalent functions which are extensively
studied in geometric function theory. Let $\mathcal{A}$ be the class of analytic functions defined on $\ID$ with the conditions $f(0)=0=f'(0)-1$ and
$\mathcal{S}$ be the class of univalent functions in $\mathcal{A}$. Thus each
$f\in \mathcal{S}$ has the following Taylor expansion about the origin
\begin{equation}\label{eq1a}
f(z)=z + \sum_{n=2}^{\infty} a_nz^n, \quad z\in \ID.
\end{equation}
Let $\mathcal{S}_k$ denote the class of functions in $\mathcal{S}$ which have
$k$-quasiconformal extension onto $\sphere$.
Finding sufficient conditions for a function $f\in \mathcal{A}$ to be univalent is a
fundamental problem in univalent function theory. We first mention here one of such conditions
proved by X. Huang (\cite{HO1}). Later, Huang and Owa (\cite[Theorem\,2]{HO}) proved similar
type of condition for quasiconformal extensions. We state these results below combining both
of them.

\begin{Thm}\label{thC}
Let $f \in \mathcal{A}$ having expansion of the form \eqref{eq1a}. If
$$
f(z)=\frac{z}{1-a_2z+\phi(z)},
$$
where $\phi(z)$ is analytic in $\ID$ with $\phi(0)=\phi'(0)=0$ and satisfies
$|(\phi(z)/z)'| \leq 1$ for all $z\in \ID$, then $f$ is univalent in $\ID$.
Moreover if $\phi$ satisfies $|(\phi(z)/z)'| \leq k$ for some
$k\in(0,1)$, then $f\in \mathcal{S}_k$ and has the following $k$-quasiconformal
extension onto  $\sphere :$
\begin{equation}\label{eq2b}
F(z)= \begin{cases}
\frac{z}{1-a_2z + \phi(z)}, & z\in \ID\\
\frac{z}{1-a_2z + |z|^2\phi(1/\overline{z})}, & z \in \IDsb.
\end{cases}
\end{equation}
\end{Thm}

Next we consider the class $\mathcal{U}_{\lambda}$ which consists of functions
$f\in \mathcal{A}$ satisfying the condition $|U_f(z)|< \lambda$,
where $\lambda \in (0,1]$ and $U_f(z):= (z/f(z))^2f'(z) - 1,~ z \in \ID$.
It is well known that $\mathcal{U}_{\lambda} \subsetneq \mathcal{S}$
(see \cite[Theorem 2]{ON}). Here we urge the readers to go through the articles
\cite{OP,OPW} and references therein for a detailed study about this class of
function. In this article we will discuss about
quasiconformal extension of each $f\in \mathcal{U}_{\lambda}$ for which
$U_f(z)\not\equiv 0$ using  Theorem C.

Now we consider the class of meromorphic functions defined on $\ID$ where the functions have pole at some nonzero point $p\in (0,1)$. Let $\mathcal{A}(p)$ denote the class of meromorphic functions $f$ that are analytic in $\ID \setminus \{p\}$, possessing a simple pole only at
$z=p$ with nonzero residue $m$ normalized by the conditions $f(0)=0=f'(0)-1$. Let
$\mathcal{S}(p):= \{ f\in \mathcal{A}(p): f $ is univalent in $\ID$\}.
Each function $f\in \mathcal{A}(p)$ has the following Taylor expansion about the
origin
\begin{equation}\label{eq2aa}
f(z)=z + \sum_{n=2}^{\infty} a_nz^n, \quad z \in \ID_p,
\end{equation}
where $\ID_p= \{z : |z|<p\}$. An important subclass of $\mathcal{S}(p)$
denoted by $\mathcal{V}_p(\lambda)$ has recently been introduced in \cite{BP}
(see also \cite{BP1}). This class can be thought as the meromorphic
analogue of the class $\mathcal{U}_{\lambda}$. We define it as below :
$$
\mathcal{V}_p(\lambda):=\{f \in \mathcal{A}(p) : |U_f(z)|< \lambda\}, \quad
\lambda\in(0,1].
$$
It is proved in \cite{BP} that $\mathcal{V}_p(\lambda)\subsetneq \mathcal{S}(p)$.
Likewise for functions in the class $\mathcal{U}_{\lambda}$, we shall also discuss quasiconformal extension of functions in $\mathcal{V}_p(\lambda)$.

Let $\mathcal{M}$ be the class of meromorphic functions defined on
$\IDs= \{z : |z|>1 \}$, such that each function
$g \in \mathcal{M}$ is analytic in $\IDs$ except for a simple pole at
$\zeta=\infty$ with $g'(\infty)=1$ and has the following expansion
\begin{equation}\label{eq2c}
g(\zeta)= \zeta + \sum_{n=0}^{\infty} b_n\zeta^{-n}, \quad \zeta \in \IDs.
\end{equation}
Let $\Sigma$ be the class of univalent functions in $\mathcal{M}$ and
$\Sigma_k$ be the class of functions in $\Sigma$ which have $k$-quasiconformal extension
onto $\sphere$.
Aksent{\' e}v\,(\cite{AK}) proved a sufficient condition for univalence for a function in
$\mathcal{M}$, which states that if $g \in \mathcal{M}$ and satisfies
$|g'(\zeta)-1|\leq 1$, for $\zeta \in \IDs$, then $g$ is univalent.
As for any $f\in \mathcal{A}$ the function $g(\zeta):=1/f(1/\zeta)\in \mathcal{M}$,
we have $g'(\zeta)-1=U_f(z)$ where $z=1/\zeta \in \ID$.
Thus the condition $|U_f(z)|< \lambda$ yields univalence of $f$ by Aksent{\' e}v's result.
In (\cite[Corollary\,2]{KR}), Krzy{\.z} proved a sufficient condition for
quasiconformal extension related to this Aksent{\' e}v's criteria. Indeed, he proved that
if $g \in \mathcal{M}$ and satisfies $|g'(\zeta)-1|\leq k/|\zeta|^2$,
$\zeta \in \IDs,k \in (0,1)$, then $g\in \Sigma_k$. Later, Becker gave a simpler
proof of the above result by constructing a suitable Loewner chain (compare \cite[Theorem\,5.3]{JB2}).
Krzy{\. z}(\cite[Theorem\,1]{KR}) also proved that if $g\in \Sigma$ of the form $g(\zeta)= \zeta + w(1/\zeta)$,
where $w$ is an analytic function in $\ID$ which satisfies the condition $|w'(z)|\leq k$, 
for some $k \in (0,1)$, then $g$ has $k$-quasiconformal extension onto $\sphere$.
As a corollary of this theorem he obtained the previous result.
In this article, we give an alternative proof of the above theorem by modifying the
Loewner chain constructed by Becker in \cite{JB2}.
 


We organize the results obtained in this article as follows. In Theorem \ref{th0}, we
first prove a result concerning quasiconformal extension of functions in the class
$\mathcal{U}_{\lambda}$.
We also give an independent proof of the above fact in
Theorem \ref{th1} by constructing a suitable Loewner chain. Next we
establish quasiconformal extensibility of functions belonging to the class
$\mathcal{V}_p(\lambda)$. This is the content of Theorem \ref{th1a}.
In Theorem \ref{th2}, we prove that
the condition $|U_g(\zeta)|< k ,\,k\in(0,1)$ is sufficient for
$k$-quasiconformal extension of a function $g \in \mathcal{M}$.
We also give a sufficient condition for a function in $\mathcal{M}$ to belong to the class
$\Sigma_k$ in Corollary \ref{cor1} and for a function in $\mathcal{A}$ to belong to the class $\mathcal{S}_k$ in Theorem \ref{th5}.
Finally, in Theorem \ref{cor2} we give an alternative proof of the
Krzy{\.z}'s result, described above by constructing a suitable Loewner chain.

\section{Main Results}
We start the section with quasiconformal extensibility of functions belonging to the class
$\mathcal{U}_{\lambda}$.
\begin{thm}\label{th0}
Each function $f\in \mathcal{U}_{\lambda}$ having expansion of the form \eqref{eq1a}
with $U_f \not\equiv 0$, where $\lambda \in (0,1)$ has $\lambda$-quasiconformal extension
onto the extended complex plane $\sphere$. The explicit extension is given by \eqref{eq2b}.
If $f \in \mathcal{U}_{\lambda}$ with
$U_f(z) \equiv 0$, then $f$ has
$|a_2|$-quasiconformal extension $(a_2\in \ID\setminus \{0\})$ onto $\sphere$ of the form
\begin{equation}\label{eq2d}
F(z)= \begin{cases}
z/(1-a_2z), & z\in \ID\\
\frac{z(|z|^2-a_{2} z)}{(|z|-a_{2} z)^2}, & z \in \IDsb .
\end{cases}
\end{equation}
\end{thm}

\bpf \underline{\textbf{Case 1\,}:} ($f\in \mathcal{U}_{\lambda}$  and $U_f(z)\not\equiv 0$)

As $f \in \mathcal{U}_{\lambda}$, we have $|U_f(z)|< \lambda$ for $z\in \ID$. Now we see
that this inequality may be replaced by $|U_f(z)| \leq \lambda |z|^2$ since
the inequalities $|U_f(z)|< 1$ and $|U_f(z)|\leq ~|z|^2$, for
$z\in \ID$ are equivalent. To prove this we observe that the analytic function
$U_f(z)$ satisfies the condition $U_f(0)=0=U_f'(0)$. Now if $|U_f(z)|<1$, then by
the Schwarz lemma, we have $|U_f(z)| \leq |z|^2$. The other way implication is trivial.
Next we observe that if $f\in \mathcal{A}$ having an expansion of the form \eqref{eq1a},
then the function $z/f(z)$ is non vanishing and analytic in $\ID$ with
$(z/f(z))_{z=0}=1$. Hence it has an expansion of the form
\begin{equation}\label{eq2d0}
z/f(z)= 1 +b_1z + b_2z^2 + \cdots, \quad z\in \ID,
\end{equation}
where $b_1=-a_2=-f''(0)/2$. Thus, each $f\in \mathcal{A}$ can be written of the form
\begin{equation}\label{eq2d1}
f(z)=z/(1-a_2z + \phi(z)), \quad z\in \ID,
\end{equation}
where $\phi$ is defined in Theorem C. Next it is easy to see that
if $f\in \mathcal{A}$, then the conditions
$|U_f(z)|\leq \lambda|z|^2$ and $|(\phi(z)/z)'| \leq \lambda$  for $z\in \ID$
are equivalent. Indeed, if $f\in \mathcal{A}$ is of the form \eqref{eq2d1}, then $\phi$
can be written in terms of $f$ as, $\phi(z)= z/f(z)+a_2z-1$. Thus
$$
|z\phi'(z)-\phi(z)|\leq \lambda |z|^2\,~ \mbox {iff}\,~
|U_f(z)| = |(z/f(z))^2f'(z)-1|\leq \lambda|z^2|.
$$
Therefore, each $f\in \mathcal{U}_{\lambda}$ can be expressed as
$ f(z)=z/(1-a_2z+\phi(z))$, such that $|(\phi(z)/z)'| \leq \lambda $, $z\in \ID$.
Thus by an application of Theorem C, we conclude that each
$f \in \mathcal{U}_{\lambda}$ for which $U_f(z)\not\equiv 0$, can be extended to a
$\lambda$-quasiconformal map onto the extended plane $\sphere$ of the form  \eqref{eq2b}.

\vspace{.2cm}

\noindent\underline{\textbf{Case 2\,}:} ($f\in \mathcal{U}_{\lambda}$  and $U_f(z)\equiv 0$)

In this case, we observe that whenever $f(z)=z/(1-a_2z + \phi(z)) \in \mathcal{U}_{\lambda}$,
then $U_f(z) \equiv 0$ if and only if $(\phi(z)/z)' \equiv 0$, and
hence $\phi(z)=cz$ for some constant~$c$.
Since $\phi(0)=0=\phi'(0)$, we must have $\phi(z) \equiv 0$.
Thus $f$ has the form as $f(z)=z/(1-a_2z)$ with $|a_2|\leq 1$.
Consequently, $f$ must be a convex univalent mapping. Now considering the following Loewner chain
$$
f(z,t)=f(z)+(e^t-1)zf'(z), \quad~(z,t)\in \ID\times[0,\infty),
$$
for convex map (compare \cite[Theorem\,1]{IH1}),
we conclude that each $f \in \mathcal{U}_{\lambda}$ for which $U_f \equiv 0$,
has $|a_2|$-quasiconformal extension ($0<|a_2|<1$) onto
$\sphere$ and the extended map is explicitly given by \eqref{eq2d}.
\epf

\noindent\textbf{Remark.} (1) If $f\in \mathcal{U}_{\lambda}$ with $U_f(z) \equiv 0$ and $a_2=0$, we are left only with the identity map on the unit disc which has a trivial quasiconformal extension onto $\sphere$ (see f.i. \cite[Example~2.6]{IH3}).
\vspace{.1cm}

\noindent (2) If $f\in \mathcal{U}_{\lambda}$ with $U_f(z) \equiv 0$ and $|a_2|=1$, then
following similar arguments provided by Hotta (\cite[Example~2.6]{IH3}), we get the quasiconformal extension of $f$ of the form
\begin{equation}\label{eq2e}
F(re^{i\theta})=
\begin{cases}
\frac{re^{i\theta}}{1-a_2 re^{i\theta}}, & r<1\\
\frac{\psi(r)e^{i\theta}}{1-a_2 \psi(r)e^{i\theta}}, & r \geq 1,
\end{cases}
\end{equation}
where $\psi : [1,\infty) \to [1,\infty)$ is continuous, injective and bi-Lipschitz
with Lipschitz constant $M>1$, such that $\psi(1)=1$ and $\psi(\infty)=\infty$.
Here we note that $F(1/a_2)= \infty$ and $F(\infty)= -1/a_2$.
A little computation yields that for $r \geq 1$, the dilatation of $F$ satisfies
$$
|\mu_F|= \left|\frac{\psi(r)-r\phi'(r)}{\psi(r)+r\phi'(r)} \right| \leq \frac{M^2-1}{M^2+1}.
$$
Thus $F$ defined in \eqref{eq2e} is a $(M^2-1)/(M^2+1)$-quasiconformal extension
of $f$ onto $\sphere$.
\vspace{.1cm}

\noindent We explain the Theorem \ref{th0} by giving an example below :
\begin{example}
Let us consider the function
$f_{\lambda}\in \mathcal{U}_{\lambda}$, $\lambda \in (0,1)$ defined by
$$
f_{\lambda}(z)= \frac{z}{1-(1+\lambda)e^{i\theta}z + \lambda e^{2i\theta}z^2}, \quad
\theta \in [0,2\pi),\, z \in \ID,
$$
which play an important role in the study of the class $\mathcal{U}_{\lambda}$
(see \cite{OPW} and the references therein). Here $a_2=(1+\lambda)e^{i\theta}$ and
$\phi(z)= \lambda e^{2i\theta}z^2$. Therefore the $\lambda$-quasiconformal extension
of $f_{\lambda}$ follows from Theorem C and is given below :
$$
F(z)=
\begin{cases}
f_{\lambda}(z), &  z \in \ID\\
\frac{z}{1-(1+\lambda)e^{i\theta}z + \lambda e^{2i\theta} z/\overline{z}},
& z\in \IDsb.
\end{cases}
$$
\end{example}


In the next theorem we give an alternative proof of quasiconformal extensibility of
a function $f \in \mathcal{U}_{\lambda}$ for which $f''(0)=0$ by constructing
suitable Loewner chain.
\begin{thm}\label{th1}
Each function $f \in \mathcal{A}$ of the form \eqref{eq1a} with $a_2=0$ that satisfies
the condition
$$
0<|U_f(z)| \leq \lambda |z|^2, \quad\text{for all}~~ z \in \ID,
$$
where $0<\lambda <1$, has $\lambda$-quasiconformal extension onto the extended plane
$\sphere$. The extended function is given explicitly as
\begin{equation}\label{eq2a}
F(z)= \begin{cases}
f(z), & z\in \ID,\\
\frac{zf(1/\overline{z})}{z-(|z|^2-1)f(1/\overline{z})}, & z \in \IDsb .
\end{cases}
\end{equation}
\end{thm}

\bpf
We consider the following function
\begin{equation}\label{eq3}
f(z,t):= \frac{zf(e^{-t}z)}{z-(e^t-e^{-t})f(e^{-t}z)}, \quad (z,t) \in \ID \times [0,\infty).
\end{equation}
Using Theorem A, we wish to show that $f(z,t)$ is a Loewner chain.
We first note that $f_0(z)=f(z,0)=f(z)$. It is easy to see that
$f(0,t)=0$ and $f'(0,t)=\lim_{z\to 0}f'(z,t)=e^{t}$ for all $t\geq 0$.
Since $f\in \mathcal{A}$ has an expansion of the form \eqref{eq1a} and by the hypothesis $a_2=0$, therefore we get	
\begin{equation}\label{eq3a}
f(z,t)= (e^tz + a_3e^{-t}z^3+ \cdots)[1-a_3(1-e^{-2t})z^2- \cdots]^{-1}= e^tz + \cdots.
\end{equation}
Hence there exists a positive constant $r_0\in (0,1)$ such that the function $f_t(z):=f(z,t)$
is analytic in $z \in \ID_{r_0}$ for all $t\geq 0$. Now it is clear that $f(z,t)$ is
absolutely continuous in $t$ for each $z \in \ID_{r_0}$. From \eqref{eq3a}, we see that
$$
\lim_{t \to \infty} e^{-t}f(z,t) = z/(1-a_3z^2),
$$
locally uniformly in $z$-variable. Hence the family
$\{e^{-t}f(z,t)\}_{t\geq 0}$ is a normal family. Therefore there exists a positive constant
$K_0$ such that $|f(z,t)|\leq K_0e^t$ for all $z\in \ID_{r_0}$ and $t \in [0, \infty)$.
Thus the first condition of Theorem A is satisfied. Next we see that a little calculation yields
\begin{equation}\label{eq4}
p(z,t):=\frac{\dot{f}(z,t)}{zf'(z,t)}= \frac{-z^2f'(e^{-t}z)e^{-t}+(e^t+e^{-t})f^2(e^{-t}z)}{z^2f'(e^{-t}z)e^{-t}-(e^t-e^{-t})f^2(e^{-t}z)}.
\end{equation}
We now show that the function $p(z,t)$ is analytic in $\ID$. Indeed, if the denominator of the last term in \eqref{eq4} vanishes for a point $z_0 \in \ID$ and $t \in[0,\infty)$,
then by the given condition we have
$$
e^{-2t}=|U_f(e^{-t}z_0)|\leq e^{-2t}|z_0|^2,
$$
which gives $|z_0|\geq 1$\,- a contradiction. Thus $p(z,t)$ is analytic in $\ID$ and $p(0,t)=1$.
Now by the hypothesis we get from \eqref{eq4} that
\begin{equation}\label{eq5}
\left|\frac{p(z,t)-1}{p(z,t)+1}\right|= e^{2t}\left|\frac{(e^{-t}z)^2f'(e^{-t}z)}{f^2(e^{-t}z)} -1 \right| \leq \lambda|z|^2< |z|^2 < 1,
\end{equation}
which implies $\text{Re}\, p(z,t)>0$. Lastly it is obvious that $p(z,t)$ is measurable
in $t$ for each fixed $z$.
From the above discussion it is clear that all the conditions in Theorem A
are satisfied and hence $f(z,t)$ defined in \eqref{eq3} is a Loewner chain.
Thus for each $t\geq 0$ the
functions $f_t$ are univalent in $\ID$ and hence in particular $f_0=f$ is univalent.
Again as $|U_f(z)| \leq \lambda |z|^2$, (which implies $f \in \mathcal{U}_{\lambda}$)
then from the equation \eqref{eq5} we get
$$
p(z,t) \in \mathcal{D}(\lambda)= \{w:|w-1|/|w+1| \leq \lambda \}.
$$
As $\mathcal{D}(\lambda)$ is contained properly in the right half plane, $f(z,t)$
is a Loewner chain which satisfies the conditions of Theorem B. Thus we conclude that each function $f(=f_0)$ in $\mathcal{U}_{\lambda}$ with $a_2=0$, (except the identity map)
has $\lambda$-quasiconformal extension to the whole plane and from \eqref{eq2} we get the extended map explicitly which is described in \eqref{eq2a}.
\epf

\noindent Next we explain the Theorem \ref{th1} through the following example, where
$f\in \mathcal{U}(\lambda)$ with $a_2:=f''(0)/2=0$.
\begin{example}
Let
$$
f(z)= z/(1 + \lambda z^2), \quad \lambda\in(0,1), ~~ z \in \ID.
$$
For this function $a_2=0$ and $U_f(z)=-\lambda z^2$ and hence
$f \in \mathcal{U}_{\lambda}$. The following $\lambda$-quasiconformal extension of $f$
can be deduced from \eqref{eq2a}:
$$
F(z)=
\begin{cases}
f(z), &  z \in \ID\\
\frac{z}{1 + \lambda z/\overline{z}},
& z\in \IDsb.
\end{cases}
$$
\end{example}

\noindent\textbf{Remark.}
(1) The $\lambda$-quasiconformal extension of the functions in
$\mathcal{U}_{\lambda}$, ($\lambda \in (0,1)$) given by \eqref{eq2a} is same as that one obtained in \eqref{eq2b}.

\vspace{.1cm}
\noindent(2) 
The \textit{Koebe function} $k(z)=z/(1-z)^2 \in \mathcal{U}(1)$, ($\lambda=1$).
We know that $k(z)$ has no quasiconformal extension to $\sphere$ (compare \cite[Example\,2.9]{IH3})). Indeed,
$k$ maps $\ID$ onto $\IC \setminus (-\infty,-1/4]$. If $k$ can be extended
quasiconformally to the whole plane, then there would exist a homeomorphism from
$\IDs$ onto $(-\infty,-1/4]$, which is not possible.

Next we obtain quasiconformal extension for functions belonging to the class
$\mathcal{V}_p(\lambda)$~:
\bthm\label{th1a}
Each function $f \in \mathcal{V}_p(\lambda)$, $0<\lambda<1$,
having the expansion of the form \eqref{eq2aa} with
$U_f(z)\not\equiv 0$, has $\lambda$-quasiconformal extension onto
$\sphere$, and is given by \eqref{eq2b}.
\ethm
\bpf
We follow the lines of proof of Theorem\,\ref{th0}.
If $f\in \mathcal{A}(p)$ with the expansion \eqref{eq2aa} then the function
$z/f(z)$ is analytic in $\ID$, non vanishing in $\ID \setminus \{p\}$ and hence it has an expansion of the form \eqref{eq2d0}. Thus each $f\in \mathcal{A}(p)$ also can be expressed
of the form \eqref{eq2d1} with $\phi(p)=a_2p-1$.
Now it is easy to see that the conditions $|U_f(z)|\leq \lambda|z|^2$
(which implies $f \in \mathcal{V}_p(\lambda)$) and
$|(\phi(z)/z)'| \leq \lambda$, for $z\in \ID$ are equivalent. Hence using
Theorem C we obtain $\lambda$-quasiconformal extension of functions in
$\mathcal{V}_p(\lambda)$ for which $U_f(z)\not\equiv 0$ and the extended function
has the form  given by \eqref{eq2b}.
\epf

\noindent\textbf{Remark.} (1) The function $f(z)=pz/(p-z)$ is the only function in
$\mathcal{V}_p(\lambda)$ for which $U_f(z)\equiv 0$ (in other words $\phi(z)\equiv 0$).
Hence by the similar arguments due to Hotta (\cite[Example\,2.6]{IH3}), we get
$(M^2-1)/(M^2+1)$-quasiconformal extension
of $f$ onto $\sphere$ as
$$
F(re^{i\theta})=
\begin{cases}
\frac{pre^{i\theta}}{p - re^{i\theta}}, & r<1\\
\frac{p \psi(r)e^{i\theta}}{p- \psi(r)e^{i\theta}}, & r \geq 1,
\end{cases}
$$
where $\psi$ and $M$ are described in \eqref{eq2e}. Here we note that $F(p)=\infty$
and $F(\infty)=-p$.

\vspace{.1cm}
\noindent(2) Let us consider the function
$$
k_p(z)= \frac{pz}{(p-z)(1-pz)}, \quad  z\in \ID.
$$
It is easy to check that $k_p \in \mathcal{V}_p(1)$. It also belongs to the class
$Co(p)$, the class of \textit{concave univalent functions with a pole} $p\in (0,1)$
(see \cite{AW,BPW}). It maps $\ID$ onto $\sphere \setminus [-p/(1-p)^2,-p/(1+p)^2]$.
Therefore, by the same reason provided by Hotta (\cite[Example\,2.9]{IH3}), $k_p$ does
not have any quasiconformal extension to the extended plane.

\vspace{.1cm}
We now explain the Theorem \ref{th1a} with the following example :
\begin{example}
Let
$$
k_p^{\lambda}(z)= \frac{pz}{(p-z)(1-\lambda pz)}
= \frac{z}{1-(\lambda p + 1/p)z + \lambda z^2}, \quad \lambda \in (0,1),~ z\in \ID.
$$
This function serves as an extremal function for many problems in the class
$\mathcal{V}_p(\lambda)$ (see f.i.\,\cite{BP}). Here $a_2=\lambda p + 1/p$ ~and
$\phi(z)=\lambda z^2$. Thus applying Theorem C we get
the $\lambda$-quasiconformal extension of $k_p^{\lambda}$ as
$$
F(z)=
\begin{cases}
k_p^{\lambda}(z), & z \in \ID \\
\frac{z}{1-(\lambda p + 1/p)z + \lambda z/\overline{z}}, & z\in \IDsb .
\end{cases}
$$
\end{example}

In Theorem \ref{th1} we have seen that the condition $|U_f(z)|< \lambda$,
$\lambda \in (0,1)$ is sufficient for univalence as well as $\lambda$-quasiconformal extensibility for a function $f\in \mathcal{A}$. In the next theorem we
prove an analogues result for functions $g\in \mathcal{M}$.

\begin{thm}\label{th2}
Let $g \in \mathcal{M}$ be such that
$$
|U_g(\zeta)| \leq k, \quad\text{for all}~\zeta \in \IDs,
$$
for some $k\in (0,1)$, then $g$ is univalent in $\IDs$ and has the following
$k$-quasiconformal extension onto $\sphere :$
\begin{equation}\label{eq7}
G(\zeta)=
\begin{cases}
g(\zeta), & \zeta \in \IDs \\
\frac{g(1/\overline{\zeta})}{1+ (1/\zeta - \overline{\zeta})g(1/\overline{\zeta})},
& \zeta \in \IDb=\{z : |z|\leq 1\}.
\end{cases}
\end{equation}
\end{thm}

\bpf
Let us consider the function
\begin{equation}\label{eq7a1}
f(z,t)= [g(e^t z\inv)]\inv + (e^t-e^{-t})z, \quad (z,t) \in \ID \times [0,\infty).
\end{equation}
Since $g\in \mathcal{M}$, a little calculation shows that $f(0,t)=0$ and $f'(0,t)=e^t$
for all $t\geq 0$. As $g\in \mathcal{M}$ has the expansion of the form \eqref{eq2c},
we have
\begin{equation}\label{eq7a}
f(z,t)= e^{-t}z[1+b_0e^{-t}z+b_1e^{-2t}z^2+\cdots]\inv + (e^t-e^{-t})z= e^tz+ \cdots.
\end{equation}
Hence there exists a positive constant $r_0\in (0,1)$ such that the function $f_t(z):=f(z,t)$
is analytic in $\ID_{r_0}$ (w.r.t. $z$-variable) for all $t\geq 0$. Again $f(z,t)$ is
absolutely continuous in $t$ for each $z \in \ID_{r_0}$. From \eqref{eq7a} we get
$$
\lim_{t \to \infty} e^{-t}f(z,t) = z,
$$
locally uniformly in $z$-variable. Hence the family $\{e^{-t}f(z,t)\}_{t\geq 0}$ is a normal family. Therefore there exists a positive constant $K_0$ such that $|f(z,t)|\leq K_0e^t$ for all $z\in \ID_{r_0}$ and $t \in [0, \infty)$.
Next a computation yields that
$$
p(z,t):=\frac{\dot{f}(z,t)}{zf'(z,t)}=\frac{-g'(e^tz\inv)e^tz\inv +
(e^t+e^{-t})zg^2(e^tz\inv)}{g'(e^tz\inv)e^tz\inv + (e^t-e^{-t})zg^2(e^tz\inv)}.
$$
Now if the denominator of the last term of the above expression vanishes for any
$z_0\in \ID$ and $t\in [0,\infty)$, then we have
$|U_g(w_0)|=e^{2t} \geq 1$ where $w_0=e^tz_0\inv \in \IDs$, which is not possible due to
our assumption. Hence $p(z,t)$ is analytic in $\ID$. Also $p(0,t)=1$ for all $t$.
Again by the given condition, we get
$$
\left|\frac{p(z,t)-1}{p(z,t)+1}\right|= e^{-2t}|U_g(w)|\leq k <1,
$$
for all $z\in \ID$ and $t\geq 0$, where $w=e^tz\inv \in \IDs$, which yields
$p(z,t)\in \mathcal{D}(k)$ (which also imply $\mathrm{Re}~p(z,t)>0$). Thus we see
that $f(z,t)$ defined in \eqref{eq7a1} satisfies all the conditions of Theorem A
and hence it becomes a Loewner chain. Thus $f(z,0)=1/g(1/z)$ is univalent in $\ID$,
which implies the univalence of $g$ in $\IDs$. Finally from Theorem B,
we have $f(z,0)=1/g(1/z)\in \mathcal{S}_k$ and hence $g(\zeta)\in \Sigma_k$,
where $\zeta=1/z\in \IDs$
and the $k$-quasiconformal extension of $g$ is given by \eqref{eq7}.
\epf

In the next corollary we establish a sufficient condition for a function in
$\mathcal{M}$ to belong to the class $\Sigma_k$.
\bcor\label{cor1}
If $g\in \mathcal{M}$ and for some $k\in(0,1)$ it satisfies
$$
\left| (\zeta /g(\zeta))^2 g'(\zeta)+1 \right| \leq k,
\quad\text{for all}~ \zeta \in \IDs,
$$
then $g\in \Sigma_k$ and has the following $k$-quasiconformal extension :
\begin{equation}\label{eq7b}
G(\zeta)=
\begin{cases}
g(\zeta), & \zeta \in \IDs \\
\frac{g(1/\overline{\zeta})}{1 - (1/\zeta - \overline{\zeta})g(1/\overline{\zeta})},
& \zeta \in \IDb.
\end{cases}
\end{equation}
\ecor
\bpf
We consider the function
$$
f(z,t)= [g(e^t z\inv)]\inv - (e^t-e^{-t})z, \quad (z,t)\in \ID \times [0,\infty).
$$
As $g\in \mathcal{M}$ with the expansion \eqref{eq2c}, we have
$$
f(z,t)= e^{-t}z[1+b_0e^{-t}z+b_1e^{-2t}z^2+\cdots]\inv - (e^t-e^{-t})z= (2e^{-t}-e^t)z+ \cdots.
$$
Thus $f(0,t)=0$ and $f'(0,t)=a_1(t):=2e^{-t}-e^t$ for all $t\geq 0$.
Therefore $|a_1(t)|\to \infty$ as $t\to \infty$.
The above expansion of $f_t(z):=f(z,t)$ also yields that there exists a number
$r_0 \in(0,1)$
such that each function $f_t$ is analytic in $z \in\ID_{r_0}$ for all $t\geq 0$.
Again we see that
$$
\lim_{t \to \infty} f(z,t)/a_1(t)=z,
$$
which implies that
$\{f(z,t)/a_1(t)\}_{t\geq 0}$ is a normal family and hence locally uniformly bounded on
$\ID$. Thus there exists a positive constant $K_0$ such that $|f(z,t)|\leq K_0|a_1(t)|$
for all $z\in \ID_{r_0}$ and $t\geq 0$. Following the similar lines of proof as given in
Theorem \ref{th2} and from the given condition, we see that the Herglotz function
$p(z,t)=\dot{f}(z,t)/zf'(z,t)$ is analytic in $z\in\ID$ and satisfies the inequality
$$
\left|\frac{p(z,t)-1}{p(z,t)+1}\right|= e^{-2t}|(w/g(w))^2g'(w)+1|\leq k <1,
$$
for all $z\in \ID$ and $t\geq 0$, where $w=e^tz\inv \in \IDs$.
Hence by \cite[Theorem~A$^{\prime}$,~Theorem~C$^{\prime}$]{IH2} we
conclude that $g\in \Sigma$ and has $k$-quasiconformal extension onto $\sphere$
of the form \eqref{eq7b}, i.e. $g \in \Sigma_k$.
\epf

\noindent\textbf{Remark.} In 1984, Brown (\cite[Theorem\,5]{JB}) proved that if
$f\in \mathcal{A}$ and satisfies $|\lambda f'(z)-1|\leq k<1,~ z \in \ID$, for some complex constant $\lambda(\neq0)$, then $f$ is univalent in $\ID$ and has $k$-quasiconformal
extension onto $\sphere$ of the form
\begin{equation}\label{eq6}
F(z)=
\begin{cases}
f(z), & z\in \ID\\
f(1/\overline{z}) + \lambda\inv (z-1/\overline{z}), & z\in \IDsb.
\end{cases}
\end{equation}
Later Hotta (see \cite[Remark\,3-1]{IH1}) proved the above result for $\lambda =1$
by constructing the Loewner chain:
$$
f(z,t)= f(e^{-t}z)+z(e^t-e^{-t}),  \quad(z,t)\in \ID\times [0,\infty).
$$
If $g\in \mathcal{M}$ such that it never vanishes on $\IDs$, then
$f(z)=1/g(1/z)\in \mathcal{A}$ for $z\in \ID$. This shows that the conditions
$|f'(z)-1|\leq k$ and $|U_g(\zeta)| \leq k~(\zeta=1/z)$ are equivalent.
Hence Theorem \ref{th2} also follows from \cite[Remark\,3-1]{IH1}
and the extended map $F$ in \eqref{eq6} (for $\lambda =1$) is related to that
one obtained in \eqref{eq7} by $F(z)=1/G(\zeta)$.

A slight modification of the Loewner chain in \cite[Remark\,3-1]{IH1} yields the following result.
\bthm\label{th5}
If $f\in \mathcal{A}$ and for some $k\in(0,1)$ it satisfies
$$
|f'(z)+ 1| \leq k, \quad\text{for all}~z\in \ID,
$$
then $f\in \mathcal{S}_k$.
\ethm
\bpf
Let us consider the function
$$
f(z,t)=f(e^{-t}z) - z(e^t-e^{-t}), \quad (z,t)\in \ID \times [0,\infty).
$$
Now it is easy to see that $f(z,0)=f(z),\,f(0,t)=0$ and $f'(0,t)=2e^{-t}-e^t=:a_1(t)$,
for all $t$. Therefore $|a_1(t)|\to \infty$ as $t\to \infty$. Since $f \in \mathcal{A}$ and
$$
\lim_{t \to \infty} f(z,t)/a_1(t)=z,
$$
therefore, $\{f(z,t)/a_1(t)\}_{t\geq 0}$ is a normal family and hence locally uniformly
bounded on $\ID$.
A straightforward calculation shows that
$$
p(z,t):=\frac{\dot{f}(z,t)}{zf'(z,t)} = \frac{-f'(e^{-t}z)e^{-t} - (e^t+e^{-t})}
{f'(e^{-t}z)e^{-t} - (e^t-e^{-t})}.
$$
If the denominator of the last term of the above expression vanishes for a point
$z_0 \in \ID$ and $t \in[0,\infty)$, then we have
$$
|f'(e^{-t}z_0)+1| = e^{2t} \geq 1,
$$
which is a contradiction. Hence $p(z,t)$ is analytic in $\ID$. Again
$$
\left|\frac{p(z,t)-1}{p(z,t)+1}\right|= e^{-2t}|f'(e^{-t}z)+1|\leq k <1,
$$
for all $z\in \ID$ and $t\geq 0$, which implies $p(z,t)\in \mathcal{D}(k)$.
Hence by \cite[Theorem~A$^{\prime}$,~Theorem~C$^{\prime}$]{IH2} we
conclude that $f\in \mathcal{S}$ has $k$-quasiconformal extension to $\sphere$, i.e.
$f\in \mathcal{S}_k$. The $k$-quasiconformal extension of $f$ has the form :
$$
F(z)=
\begin{cases}
f(z), & z \in \ID \\
f(1/\overline{z})-z + 1/\overline{z}, & z\in \IDsb.
\end{cases}
$$
\epf


As mentioned in the previous section we will see in the next theorem that the
univalence and quasiconformal extensibility criteria due to Krzy{\.z}
can be proved by constructing suitable Loewner chain.

\bthm\label{cor2}
Let $g \in \mathcal{M}$ of the form $g(\zeta)=\zeta+ w(1/\zeta)$ for $\zeta\in \IDs$,
where $w$ is an analytic function in $\ID$ with $w(0)=0$ and satisfies
$|w'(z)|\leq 1$ for all $z\in \ID$. Then $g$ is univalent in $\IDs$. Moreover if
$w$ satisfies $|w'(z)|\leq k$, for some $k \in (0,1)$ then $g$ has $k$-quasiconformal
extension onto $\sphere$, that is given by
$$
G(\zeta)=
\begin{cases}
g(\zeta) = \zeta + w(1/\zeta), & \zeta \in \IDs \\
\zeta + w(\overline{\zeta}), & \zeta \in \IDb.
\end{cases}
$$
\ethm
\bpf
In \cite{JB2}, Becker considered the following Loewner chain 
\begin{equation*}
f(z,t)= [g(e^tz\inv)-(e^t-e^{-t})z\inv]\inv, \quad (z,t) \in \ID \times [0,\infty).
\end{equation*}
Now putting $g(\zeta)=\zeta+ w(1/\zeta)$, we get from above
\begin{equation}\label{eq9}
f(z,t)= [w(e^{-t}z) + e^{-t}z\inv]\inv,\quad (z,t)\in \ID \times [0,\infty).
\end{equation}
We may assume $b_0=0$ in the expansion \eqref{eq2c} for $g\in \mathcal{M}$
which is equivalent to the condition $w(0)=0$, then $w$ has an expansion of the form
$$
w(z)= b_1z + b_2z^2 + b_3z^3 + \cdots, \quad\text{for} ~z\in \ID.
$$
Now it is easy to show that the function defined in \eqref{eq9} is a Loewner chain.
Indeed, we see that $f(0,t)=0$ and $f'(0,t)=e^t$ for all $t\geq 0$ and
$$
f(z,t)= e^tz[1+b_1z^2+b_2e^{-t}z^3+ \cdots]\inv = e^tz + \cdots,
$$
that implies
$$
\lim_{t \to \infty} e^{-t}f(z,t) = z/(1+b_1z^2),
$$
locally uniformly in $z\in \ID$.

Thus following similar kind of arguments as provided in Theorem \ref{th2},
$f(z,t)$ in \eqref{eq9} turns out to be a Loewner chain from Theorem A. 
Now the quasiconformal extension of
$f(z,0)=1/g(1/z)$ and hence of $g(\zeta)$, $(\zeta=1/z)$ follows from the Theorem B.
\epf

\noindent\textbf{Remark.}
The condition $w(0)=0$ in Theorem\,\ref{cor2} is necessary for $f(z,t)$ defined in
\eqref{eq9} to be a Loewner chain. This condition was not needed in \cite[Theorem\,1]{KR}.

\end{document}